\documentclass{amsart}
\usepackage{amsmath, amssymb}
\usepackage{xcolor}

\input xy
\xyoption{all}

\newtheorem{innercustomthm}{Corollary}
\newenvironment{customthm}[1]
  {\renewcommand\theinnercustomthm{#1}\innercustomthm}
  {\endinnercustomthm}

\newtheorem{theorem}{Theorem}[section]
\newtheorem*{theorem*}{Theorem}
\newtheorem{lemma}[theorem]{Lemma}
\newtheorem{corollary}[theorem]{Corollary}
\newtheorem*{corollary*}{Corollary}

\newtheorem{proposition}[theorem]{Proposition}
\newtheorem{claim}[theorem]{Claim}

\theoremstyle{definition}

\newtheorem{remark}[theorem]{Remark}
\newtheorem{definition}[theorem]{Definition}

\def\Frac{\operatorname{Frac}}
\def\alg{\operatorname{alg}}

\def\id{\operatorname{id}}
\def\span{\operatorname{span}}

\def\fix{\operatorname{Fix}}
\def\adim{\operatorname{adim}}
\def\eadim{\operatorname{sadim}}
\def\tdim{\operatorname{tdim}}
\def\aut{\operatorname{Aut}}
\def\bir{\operatorname{Bir}}
\def\gal{\operatorname{Gal}}
\def\hilb{\operatorname{Hilb}}

\def\AA{\mathbb{A}}

\def\NN{\mathbb{N}}
\def\KK{\mathbb{K}}

\def\G{\mathcal{G}}

\def\dto{\dashrightarrow}

\begin{document}

\title[On invariant rational functions under rational transformations]{On invariant rational functions under rational transformations}

\author{Jason Bell}
\address{Jason Bell\\
University of Waterloo\\
Department of Pure Mathematics\\
200 University Avenue West\\
Waterloo, Ontario \  N2L 3G1\\
Canada}
\email{jpbell@uwaterloo.ca}

\author{Rahim Moosa}
\address{Rahim Moosa\\
University of Waterloo\\
Department of Pure Mathematics\\
200 University Avenue West\\
Waterloo, Ontario \  N2L 3G1\\
Canada}
\email{rmoosa@uwaterloo.ca}

\author[Matthew Satriano]{Matthew Satriano}
\address{Matthew Satriano\\
University of Waterloo\\
Department of Pure Mathematics\\
200 University Avenue West\\
Waterloo, Ontario \  N2L 3G1\\
Canada}
\email{msatrian@uwaterloo.ca}

\date{\today}

\thanks{The authors were supported by Discovery Grants from the National Sciences and Engineering Research Council of Canada.
R. Moosa and M. Satriano were also supported by Waterloo Math Faculty Research Chairs.}



\begin{abstract}
Let $X$ be an algebraic variety equipped with a dominant rational map $\phi:X\dashrightarrow X$.
A new quantity measuring the interaction of $(X,\phi)$ with trivial dynamical systems is introduced; 
the {\em stabilised algebraic dimension} of $(X,\phi)$ captures the maximum number of new algebraically independent invariant rational functions on $(X\times Y,\phi\times\psi)$, as $\psi:Y\dashrightarrow Y$ ranges over all dominant rational maps on algebraic varieties.
It is shown that this birational invariant agrees with the maximum $\dim X'$ where $(X,\phi)\dashrightarrow(X',\phi')$ is a dominant rational equivariant map and $\phi'$ is part of an algebraic group action on $X'$.
As a consequence, it is deduced that if some cartesian power of $(X,\phi)$ admits a nonconstant invariant rational function, then already the second cartesian power does.
\end{abstract}

\maketitle

%\tableofcontents

\section{Introduction}
\noindent
Given an irreducible quasi-projective variety~$X$ over a field~$k$, by a {\em rational transformation} of $X$ we will mean a dominant rational map $\phi:X\dashrightarrow X$ over $k$, and we will refer to the pair $(X,\phi)$ as a {\em rational dynamical system}.
A rational transformation induces, and is determined by, a difference field structure on the function field, namely the $k$-endomorphism  $\phi^*:k(X)\to k(X)$ given by $\phi^*(f)=f\circ \phi$.
An {\em invariant rational function} is an element of the fixed field $\fix(k(X),\phi^*)$.

%A natural first step in the birational classification problem, would be to measure how far $(X,\phi)$ is from being a trivial dynamical system:

The following is a fundamental invariant in the study of rational dynamical systems that measures how far $(X,\phi)$ is from being trivial (or rather periodic):

\begin{definition}
\label{defn:ad}
The {\em algebraic dimension} of $(X,\phi)$, denoted by $\adim(X,\phi)$, is the transcendence degree of $\fix(k(X),\phi^*)$ over $k$.
That is, it is the maximum number of algebraically independent invariant rational functions.
\end{definition}

Our terminology here is borrowed from the theory of compact complex manifolds, with which there is a loose analogy.

Although this invariant has not been previously named,  it appears throughout the dynamics literature.
The condition of having algebraic dimension zero is especially well-studied, for example in~\cite{brs} and~\cite{Can} where it is shown to imply (in characteristic zero) the existence of only finitely many invariant hypersurfaces.
An important consequence for us of $\adim(X,\phi)=0$, in the case that $k=\mathbb C$, is a theorem of Amerik and Campana~\cite{amerik-campana} asserting that there must exist a $k$-point whose orbit under $\phi$ is Zariski dense in~$X$.
This was extended to uncountable algebraically closed fields of any characteristic in~\cite{BGR}.
(When $k=\mathbb Q^{\alg}$ this is the still open Zariski dense orbit conjecture, see~\cite{zhang, amerik, alicetom}.)

Algebraic dimension can be expressed as the maximum $\dim Y$ such that there is a dominant equivariant rational map from $(X,\phi)$ to the trivial dynamical system $(Y,\id)$.
Replacing $(X,\phi)$ by the generic fibre of such a maximal image, one can sometimes reduce the birational classification problem to the case of algebraic dimension~$0$; see Lemma~\ref{are} below.

\begin{remark}
While we will work throughout in arbitrary characteristic, the reader should be warned that the above definition of algebraic dimension is maybe too naive in positive charactetistic:
if $X$ is defined over a finite field $\mathbb F_q$, and $\phi$ is the $q$-power Frobenius, then every subvariety of $X$ over $\mathbb F_q$ is $\phi$-invariant.
We might therefore be lead to see $(X,\phi)$ as being close to the trivial dynamics $(X,\id)$, but in fact $\adim(X,\phi)=0$.
\end{remark}

A drawback of algebraic dimension is that it may not detect interactions with trivial dynamics that only appear after base extension.
For example, suppose $\operatorname{char}(k)=0$ and consider $X=\mathbb A^1$ equipped with $\phi(x)=x+1$.
Then there are no nonconstant invariant rational functions, i.e., $\adim(X,\phi)=0$.
However, after base changing by $(X,\phi)$ itself, we obtain a dynamical system on $\mathbb A^2$ given by $(x,y)\mapsto (x+1,y+1)$, which acquires the invariant rational function $x-y$.

To deal with such cases, we introduce a ``stable limit'' of algebraic dimension, which is a more robust measure:

\begin{definition}
\label{ead}
Suppose $X$ is geometrically irreducible.
The {\em stabilised algebraic dimension} of $(X,\phi)$, denoted by  $\eadim(X,\phi)$, is the maximum of
$$\adim(X\times Y,\phi\times\psi)-\adim(Y,\psi),$$
as $(Y,\psi)$ ranges over all irreducible rational dynamical systems over~$k$.
\end{definition}

Here we have assumed geometric irreducibility for convenience, so that the product $X\times Y$ is again irreducible and we can speak of its rational function field.

Thus, $\eadim(X,\phi)$ measures the maximum number of new invariant rational functions one may obtain after base change. It is clear that
\[
\adim(X,\phi)\leq \eadim(X,\phi)
\]
since one may take $(Y,\psi)$ to be the one-point system. It turns out that $\eadim(X,\phi)$ is a bounded quantity:~we show that
\[
\eadim(X,\phi)\leq \dim X
\]
in Corollary~\ref{dimxbound} below.

The two extreme cases, when $\eadim(X,\phi)=0$ and when $\eadim(X,\phi)=\dim X$, are of particular interest and have arisen before in other guises.
They correspond, roughly, to what Chatzidakis and Hrushovski call ``fixed-field-free" and ``fixed-field-internal", respectively, in~\cite[$\S1.2$]{ChHr-AD1}.

Our main theorem relates $\eadim(X,\phi)$ to a seemingly unrelated quantity arising from algebraic group actions.
To describe this latter quantity, we introduce some convenient terminology for rational transformations that come from an algebraic group action.

\begin{definition}
\label{defn:tv}
A rational transformation $\phi:X\dto X$ is called a {\em translation} if there exists an algebraic group $G$ and a faithful algebraic group action
$$\rho:G\times X\to X,$$ over $k$, such that 
$\phi=\rho(g,-)$ for some $g\in G(k)$.
In this case we say that the rational dynamical system $(X,\phi)$ is {\em translational}.\footnote{Such $(X,\phi)$ were called {\em translation varieties} in~\cite{ChHr-AD1}.} 
\end{definition}

The rational dynamical system $(\mathbb A^1,x\mapsto x+1)$ considered above is translational with $G$ equal to the additive group. 
It turns out that this simple example is typical of the general case:
if $(X,\phi)$ is translational then $\eadim(X,\phi)=\dim X$.
If, in addition, $\adim(X,\phi)=0$, then the stabilised algebraic dimension is already witnessed by considering the cartesian square $(X\times X,\phi\times\phi)$.
See Propositions~\ref{prop:sdim=dim} and~\ref{translation}, below, for these properties of translations.

The following measures how far a rational dynamical system is from being translational:

\begin{definition}
\label{defn:td}
The {\em translation dimension} of $(X,\phi)$, denoted by $\tdim(X,\phi)$, is the maximum $n\geq 0$ such that there is a translational dynamical system $(Y,\psi)$ over~$k$ of dimension~$n$, and a dominant rational equivariant map $(X,\phi)\dashrightarrow (Y,\psi)$.
\end{definition}

Our main theorem appears as Theorem~\ref{thm:tve} below:

\begin{theorem*}
Suppose $(X,\phi)$ is a rational dynamical system over an algebraically closed field $k$.
Then
$\eadim(X,\phi)=\tdim(X,\phi)$.
\end{theorem*}

An immediate consequence is that if stabilised algebraic dimension is maximal (i.e., equal to dimension) then $(X,\phi)$ is a generically finite cover of a translational dynamical system.
But translational dynamical systems are preserved under finite covers (this is Proposition~\ref{covers} below), and hence we obtain:

\begin{customthm}{A}
\label{cor:A}
The stabilised algebraic dimension of $(X,\phi)$ is equal to $\dim X$ if and only if $(X,\phi)$ is birationally equivalent to a translational dynamical system.
\end{customthm}

This appears as Corollary~\ref{maxsdim} below.

We also get an interesting consequence when the stabilised algebraic dimension is only assumed to be positive; the following appears as Corollary~\ref{orthbound} below:

\begin{customthm}{B}
\label{cor:B}
The following are equivalent:
\begin{itemize}
\item[(i)]
There exists $(Y,\psi)$ with
$\adim(X\times Y,\phi\times\psi)-\adim(Y,\psi)>0$.
\item[(ii)]
$\adim(X\times X,\phi\times\phi)>0$.
\end{itemize}
In particular, some cartesian power of $(X,\phi)$ admits a nonconstant invariant rational function if and only if the second power does.
\end{customthm}

\begin{remark}
This statement does not hold for nonabelian group actions.
For instance, since $\operatorname{PGL}_2$ acts $2$-transitively on the projective line, $\mathbb P^1$, every $\operatorname{PGL}_2$-invariant rational function on $\mathbb P^1\times\mathbb P^1$ is constant; but the cross ratio determines an invariant, nonconstant, rational function on $(\mathbb P^1)^4$.
\end{remark}

We note that Corollary \ref{cor:B} is the precise analogue in algebraic dynamics of a recent result of the second-named author and R\'emi Jaoui about  algebraic vector fields~\cite[Theorem~1.1]{abred}.
Our method of proof, however, is {\em not} analogous to \cite{abred}. 
The main tool in~\cite{abred} is the model theory of differentially closed fields, and in particular the manifestation there of internality and the corresponding binding group actions.
Algebraic dynamics may also be viewed in a model-theoretic setting; the model companion of difference fields was developed by Chatzidakis and Hrushovski in~\cite{acfa} and applied to algebraic dynamics in~\cite{ChHr-AD1} and~\cite{ChHr-AD2}, for example.
However the situation is rather more subtle, due especially to the failure of quantifier elimination in this setting.
There is a theory of binding group actions arising from quantifier-free internality, especially in the case of linear difference equations, introduced by Kamensky in~\cite{Moshe}.
But the machinery is not sufficiently developed outside the linear case to facilitate a direct adaptation of the methods of~\cite{abred}.

Our approach here, therefore, is not model-theoretic, but rather algebraic and geometric.
The key step in our proof of the Theorem is to construct, given $(X,\phi)$ of stabilised algebraic dimension $\geq n$, an $n$-dimensional rational image $(Y,\psi)$ with the property that all the iterates of $\psi$ live in a fixed algebraic family of rational transformations.
Such uniform definability of the iterates of $\psi$ forces it to be a translation (see~\cite[$\S$2]{cantat14} or Proposition~\ref{findtv} below).

The plan of the paper is as follows.
In Section~\ref{sect:ad} we study algebraic dimension more closely, showing in particular that it is at least $\dim X-n$ if there exists a covering family of $n$-dimensional invariant subvarieties (Proposition~\ref{families}).
In Section~\ref{sect:td} we study translations more closely, proving that they have maximal stabilised algebraic dimension, and establishing the above mentioned criteria for being translational.
Finally, in Section~\ref{sect:ftv} we prove the main Theorem and its Corollaries.

\medskip
{\bf Acknowledgements.}
We would like to thank the anonymous referee,
whose comments led to significant improvements in the paper, expanding our results to 
arbitrary characteristic and to dominant rational maps (rather than birational transformations).

\bigskip
\section{Algebraic dimension}
\label{sect:ad}

\noindent
We fix throughout a base field $k$.

\medskip
\subsection{Algebraic reduction}
Given rational dynamical systems $(X,\phi)$ and $(Y,\psi)$ over~$k$, and a dominant rational map $\alpha:X\dashrightarrow Y$, we say that $\alpha$ is {\em equivariant} if $\alpha\circ\phi=\psi\circ\alpha$.
This is often indicated by writing $\alpha:(X,\phi)\dashrightarrow (Y,\psi)$.

By an {\em invariant rational function} of $(X,\phi)$ we mean a rational function $f\in k(X)$ such that $\phi^*(f)=f$.
That is, $f$ is in $\fix(k(X),\phi^*)$, the fixed field of the induced difference field $(k(X),\phi^*)$.
Equivalently, $f:(X,\phi)\dashrightarrow (\mathbb A_k^1,\id)$ is equivariant.

More generally, by an {\em invariant fibration} of $(X,\phi)$ we mean an equivariant dominant rational map $(X,\phi)\dashrightarrow(Y,\id)$, for some irreducible algebraic variety $Y$ over~$k$.
By the {\em dimension} of the invariant fibration we mean the dimension of $Y$.
So a nonconstant invariant rational function is a $1$-dimensional invariant fibration.

Suppose $f:X\dashrightarrow Y$ is an invariant fibration of $(X,\phi)$, and let $F=k(Y)$.
Then the generic fibre of $f$, which we will denote by $X_f$, is a closed irreducible subvariety of $X_F$ passing through a generic point of~$X$ over $k$.
It follows that $X_f$ is not contained in the indeterminacy locus of $\phi_F$, is $\phi_F$-invariant, and the restriction $\phi_F|_{X_f}:X_f\dashrightarrow X_f$ is a dominant rational map over~$F$.
We will, hopefully without confusion, continue to denote $\phi_F|_{X_f}$ by $\phi$.
Hence $(X_f,\phi)$ is itself a rational dynamical system over $F$.

\begin{lemma}
\label{are}
Every rational dynamical system $(X,\phi)$ admits a maximal invariant fibration $f:X\dashrightarrow Y$ in the sense that all invariant fibrations factor through $f$.
This invariant fibration has the properties that its generic fibre, $(X_f,\phi)$ over $k(Y)$, has no nonconstant invariant rational functions, and that $k(Y)=\fix(k(X),\phi^*)$.
\end{lemma}

\begin{proof}
Let $F=\fix(k(X),\phi^*)$.
Then $F=k(Y)$ for some irreducible algebraic variety $Y$ over $k$ and we can consider $(Y,\id)$.
The inclusion $k(Y)\subseteq k(X)$ induces a dominant rational map $f:X\dashrightarrow Y$, which is an invariant fibration because $\phi^*$ acts trivially on $F$.
So the generic fibre of $f$ gives rise to a rational dynamical system $(X_f,\phi)$ over $F$.
Now $F(X_f)=k(X)$ so that $\fix(F(X_f),\phi^*)=F$.
That is, $(X_f,\phi)$ has no nonconstant invariant rational functions.

To see that $f$ is maximal, note that if $g:X\dashrightarrow Z$ were another invariant fibration then $k(Z)\subseteq\fix(k(X),\phi^*)=k(Y)$.
\end{proof}

The invariant fibration produced above is unique up to birational equivalence: if we have two maximal invariant fibrations, $f_i:(X,\phi)\dashrightarrow(Y_i,\id)$ for $i=1,2$, then there is a unique birational map $g:Y_1\dashrightarrow Y_2$ such that $gf_1=f_2$.

\begin{definition}
We call the maximal invariant fibration given by Lemma~\ref{are} the {\em algebraic reduction} of $(X,\phi)$, and denote it by $r_\phi$.
\end{definition}

Recall from Definition~\ref{defn:ad} that the {\em algebraic dimension of $(X,\phi)$}, denoted by $\adim(X,\phi)$, is the transcendence degree of $\fix(k(X),\phi^*)$ over $k$.
We see from the above construction that this is equal to the dimension of the algebraic reduction, that is, the dimension of the image of $r_\phi$.

When $X$ is geometrically irreducible, $\adim(X,\phi)>0$ if and only if $(X,\phi)$ has a nonconstant invariant rational function.
This is because $k$ is relatively algebraically closed in $k(X)$, so $\fix(k(X),\phi^*)$ has positive transcendence degree over $k$ if and only if it is not equal to $k$.
However, even when $X$ is geometrically irreducible it is not necessarily the case that the generic fibre of the algebraic reduction will be geometrically irreducible.
The obstruction, however, can be removed by replacing $\phi$ with some iterate:

\begin{lemma}
\label{iterate}
Suppose $(X,\phi)$ is a rational dynamical system over $k$.
There exists $N\geq 1$ such that for all positive multiples $m$ of $N$, the generic fibre of $r_{\phi^m}$ is geometrically irreducible.
\end{lemma}

\begin{proof}
Let $F=\fix(k(X),\phi^*)$ and $L=F^{\alg}\cap k(X)$.
Then $L$, being a finitely generated algebraic extension of $F$ is a finite extension of $F$, say of degree $d\geq 1$.
Let $N:=d!$.

Given any positive integer multiple $m$ of $N$, we claim that $\fix(k(X),(\phi^*)^m)=L$.
Indeed, it is clear that
$\fix(k(X),(\phi^*)^m)\subseteq\fix(k(X),\phi^*)^{\alg}\cap k(X)=L$.
For the opposite containment, note that $\phi^*$ restricts to an automorphism of $L$ over $F$, as it is an $F$-endomorphism of $L$ and $L$ is a finite extension of $F$.
Since $\aut(L/F)$ has order dividing $N$, we have that $(\phi^*)^N$ is the identity on $L$.

But $(\phi^*)^m=(\phi^m)^*$.
So if we let $Y$ be an irreducible variety over~$k$ such that $k(Y)=L$, then the induced dominant rational map $X\dashrightarrow Y$ is the algebraic reduction $r_{\phi^m}$ of $(X,\phi^m)$.
Since $k(Y)$ is relatively algebraically closed in $k(X)$, the generic fibre of~$r_{\phi^m}$ is geometrically irreducible.
\end{proof}

It may be worth recording the observation that replacing $\phi$ with an iterate does not change the algebraic dimension:

\begin{lemma}
\label{iterateOK}
Suppose $(X,\phi)$ is a rational dynamical system over $k$ and $m$ is a positive integer.
Then $\fix\big(k(X),(\phi^m)^*\big)$ is algebraic over $\fix(k(X),\phi^*)$, and hence
$\adim(X,\phi)=\adim(X,\phi^m)$.

Moreover, if the generic fibre of $r_\phi$ is geometrically irreducible then $r_\phi=r_{\phi^m}$.
\end{lemma}

\begin{proof}
Since $(\phi^m)^*=(\phi^*)^m$,
any $f\in \fix\big(k(X),(\phi^m)^*\big)$ is a root of the polynomial
$$(x-f)(x-\phi^*f)(x-(\phi^*)^2f)\cdots(x-(\phi^*)^{m-1}f)$$
 whose coefficients are in $\fix(k(X),\phi^*)$.
 Hence, $\fix\big(k(X),(\phi^m)^*\big)$ is algebraic over $\fix(k(X),\phi^*)$, and so have the same transcendence degree over~$k$.
 
 For the ``moreover" clause, note that the assumption that $X_{r_\phi}$ is geometrically irreducible means precisely that $\fix(k(X),\phi^*)$ is relatively algebraically closed in $k(X)$.
 By what we have just seen, it follows that $\fix(k(X),(\phi^m)^*)=\fix(k(X),\phi^*)$, and hence $r_\phi=r_{\phi^m}$.
\end{proof}

\medskip
\subsection{Families of invariant subvarieties}
There is a strong connection between invariant fibrations and families of invariant subvarieties.
First of all, by an {\em invariant} subvariety of a rational dynamical system $(X,\phi)$ over $k$, we mean a (closed) subvariety $Z\subseteq X$ over $k$ that is not contained in the indeterminacy locus of $\phi$ and such that $\phi(Z)\subseteq Z$.
We are interested in algebraic families of such:

\begin{definition}
Suppose $(X,\phi)$ is a rational dynamical system over $k$.
By a {\em family of $n$-dimensional invariant subvarieties} we mean an irreducible variety $V$ over $k$ and an irreducible subvariety $Z\subseteq V\times X$ such that the co-ordinate projection $Z\to V$ is dominant,
$n=\dim Z-\dim V$,
and $Z$ is invariant for $(V\times X,\id\times\phi)$.
We say the family is {\em covering} if the co-ordinate projection $Z\to X$ is also dominant.
\end{definition}

Note that we do not require $Z\to V$ to be flat.

For example, an invariant fibration $f:X\dashrightarrow Y$, or rather its graph viewed as a subvariety of $Y\times X$, is an example of a covering family of  $n$-dimensional invariant subvarieties, where $n$ is $\dim X-\dim Y$.

The terminology might require some justification.
Suppose $Z\subseteq V\times X$ is a covering family of $n$-dimensional invariant subvarieties, and let $L=k(V)$.
Then the generic fibre of $Z\to V$, which we will denote by $Z_v$, is an irreducible $n$-dimensional subvariety of $X_L$ which is $\phi_L$-invariant, and which passes through a generic point of $X$ over $k$.

Here is our main tool for recognising when algebraic dimension is high:

\begin{proposition}
\label{families}
If $(X,\phi)$ admits a covering family of $n$-dimensional invariant subvarieties then
$\adim(X,\phi)\geq \dim X-n$.
\end{proposition}

\begin{proof}
Fix $Z\subseteq V\times X$, a covering family of $n$-dimensional invariant subvarieties.
Let $r_\phi:X\dashrightarrow Y$ be the algebraic reduction of $(X,\phi)$, with generic fibre $(X_{r_\phi},\phi)$.
We will show that $n\geq \dim X_{r_\phi}$.
This is what is required as $\dim X_{r_\phi}=\dim X-\dim Y$ and $\dim Y=\adim(X,\phi)$.

Note that $Z$ is also a family of invariant subvarieties of $(X,\phi^m)$ for any $m\geq 1$.
So, by Lemma~\ref{iterateOK}, we may freely replace $\phi$ with an iterate in this proof.
In particular, by Lemma~\ref{iterate}, we may assume that $(X_{r_\phi},\phi)$ is geometrically irreducible.

Fix an uncountable algebraically closed field $\mathbb K$ that extends both $F:=k(Y)$ and $L:=k(V)$.
As $X_{r_\phi}$ is geometrically irreducible the base extension to $\mathbb K$ is still irreducible, and we have a rational dynamical system over~$\mathbb K$, which we will also denote by $(X_{r_\phi},\phi)$, hopefully without confusion.
Let $\sigma$ be the induced $\mathbb K$-endomorphism of $\mathbb K(X_{r_\phi})$.
It extends $\phi^*$ on $F(X_{r_\phi})=k(X)$.
Note that $\fix(\mathbb K(X_{r_\phi}),\sigma)=\mathbb K$ since $\fix(k(X),\phi^*)=F$.
In other words, $(X_{r_\phi},\phi)$ continues to have no nonconstant invariant rational functions over~$\mathbb K$.
Now, a rational transformation over an uncountable algebraically closed field, if it has no nonconstant invariant rational functions, must have a Zariski dense orbit -- this is \cite[Th\'eor\`eme~4.1]{amerik-campana} in the complex case and~\cite[Corollary~6.1]{BGR} for arbitrary uncountable algebraically closed fields.
So, there is $a\in X_{r_\phi}(\mathbb K)$ generic over $F$ such that the orbit of $a$ under the action of $\phi$ is Zariski dense in $X_{r_\phi}$.
Since $a$ is generic in $X$ over $k$, we have that a generic fibre $Z_v$ of $Z\to V$ passes through $a$.
Since $Z_v$ is $\phi$-invariant it contains the $\phi$-orbit of~$a$, and hence has dense intersection with $X_{r_\phi}$.
So $X_{r_\phi}(\mathbb K)\subseteq Z_v(\mathbb K)$.
It follows that
$$n=\dim Z-\dim V=\dim Z_v\geq\dim X_{r_\phi}$$
as desired.
\end{proof}

\bigskip
\section{Translations}
\label{sect:td}

\noindent
We now consider rational transformations coming from actions of algebraic groups.
Recall from Definition~\ref{defn:tv} that $\phi:X\dto X$ is a {\em translation} over a field $k$ if there is an algebraic group $G$ over $k$ with a faithful algebraic group action $G\times X\to X$ over $k$, such that $\phi$ agrees with the action of some $k$-point of $G$.

\subsection{Stabilised algebraic dimension for translations}
The stabilised algebraic dimension of a translation is maximal:

\begin{proposition}
\label{prop:sdim=dim}
If $(X,\phi)$ is a geometrically irreducible translational dynamical system then $\eadim(X,\phi)\geq\dim X$.
\end{proposition}

\begin{proof}
Let $G$ be an algebraic group acting faithfully on $X$ such that $\phi$ agrees with the action of some $g\in G(k)$.
We may assume that $G$ is the Zariski closure of the cyclic subgroup generated by~$g$.
It follows that $\fix(k(X),\phi^*)$ is equal to the field of $G$-invariants of $k(X)$.
In particular, $\adim(X,\phi)$ is equal to the transcendence degree of the field of $G$-invariants of $k(X)$ over~$k$.

Now, a theorem of Rosenlicht (see~\cite[Theorem~2]{Ros}, but also~\cite{BGR}) says that there is a $G$-invariant nonempty Zariski open subset $U\subseteq X$ on which the orbits of $G$ arise as the fibres of a dominant (separable) morphism $\tau:U\to Y$.
Moreover, $\dim Y$ is the transcendence degree of the field of $G$-invariants in $k(X)$.
Shrinking~$U$ further, if necessary, we may assume that the orbits of~$G$ in~$U$ are all of the same dimension, say $\ell$.
Hence, $\adim(X,\phi)=\dim X-\ell$.

Consider, for any cartesian power $X^n$, the diagonal action of $G$ and $\phi$ on $X^n$.
Again $\fix(k(X^n),\phi^*)$ is equal to the field of $G$-invariants of $k(X^n)$.
And again there is a $G$-invariant nonempty Zariski open subset $U_n\subseteq X^n$ on which the orbits of $G$ are all of dimension
$\ell_n=n\dim X-\adim(X^n,\phi)$.

Now, by the descending chain condition on algebraic subgroups, the action of~$G$ on~$X$ has a finite {\em base}, namely there are $n\geq 1$ and $x_1,\ldots ,x_n\in X$ such that $\{g\in G\colon g\cdot x_i = x_i~{\rm for}~i=1,\ldots ,n\}$ is $\{1\}$.
It follows that the orbit of 
$(x_1,\ldots ,x_n)$ in $X^{n}$ under the action of $G$ has dimension $\dim G$.
Moreover, one can find a {\em generic} base, namely, we may assume that $(x_1,\dots,x_n)\in U_n$.
It follows that $\ell_n=\dim G$.
As every superset of a base is again a base, the same argument gives that $\ell_m=\dim G$ for all $m\geq n$.
In particular, 
\begin{eqnarray*}
\adim(X^{n+1},\phi)-\adim(X^n,\phi)
&=&
(n+1)\dim X-\dim G-(n\dim X-\dim G)\\
&=&
\dim X.
\end{eqnarray*}
This witnesses that $\eadim(X,\phi)\geq\dim X$.
\end{proof}

In fact, we will see in Corollary~\ref{dimxbound}, below, that $\dim X$ is always an upper bound for $\eadim(X,\phi)$, and hence we have equality in the above proposition.

One of the main results of this paper, namely Corollary~\ref{cor:A} of the Introduction, which is also Corollary~\ref{maxsdim} below, is a converse to Proposition~\ref{prop:sdim=dim}.
Namely, the only way stabilised algebraic dimension can be maximal is if $(X,\phi)$ is translational.

When $\adim(X,\phi)=0$, using our Proposition~\ref{families} instead of Rosenlicht's theorem, we can do a little better than Proposition~\ref{prop:sdim=dim}; we can show that the stabilised algebraic dimension of a translation is witnessed by the second cartesian power.
This will turn out to be useful, and is a consequence of the following more general, relative, statement:

\begin{proposition}
\label{translation}
Suppose $(X,\phi)$ is translational and the generic fibre, $X_{r_\phi}$, of the algebraic reduction, $r_\phi:X\dashrightarrow Y$, is geometrically irreducible.
Then
$$\adim(X_{r_\phi}\times X_{r_\phi},\phi\times\phi)\geq\dim(X_{r_\phi}).$$
\end{proposition}

\begin{proof}
The assumption that $X_{r_\phi}$ is geometrically irreducible ensures that the product $X_{r_\phi}\times X_{r_\phi}$ is again irreducible, whose algebraic dimension under $\phi\times\phi$ it makes sense to compute.

Fix $m\geq 1$.
Note that $(X,\phi^m)$ is again translational.
By Lemma~\ref{iterateOK}, we have that $r_\phi=r_{\phi^m}$, and
$\adim(X_{r_{\phi^m}}\times X_{r_{\phi^m}},\phi^m\times\phi^m)=\adim(X_{r_{\phi}}\times X_{r_{\phi}},\phi\times\phi)$.
So, in order to prove the proposition, we may freely replace~$\phi$ with an iterate.

We are given that $\phi$ acts on $X$ as multiplication by $g\in G(k)$ for some algebraic group $G$ over $k$ acting faithfully on $X$.
Let $H$ be the Zariski closure of the subgroup of $G$ generated by $g$; it is a commutative  algebraic subgroup of $G$ over~$k$.
Replacing $\phi$ with an iterate again if necessary, and hence $g$ by a power of~$g$, we may assume that $H$ is connected.
Consider the action of $H$ on $X$, as well as its base extension to $F:=k(Y)$.
Note that  $X_{r_\phi}$ is $H$-invariant.

%Indeed, Let $\mathbb K$ be an uncountable algebraically closed field extending $F$, for any $x\in X(\mathbb K)$, we see that$$\{h\in G:hx\in X_{r_\phi}\}$$ is a Zariski closed subset over $F(x)$. Hence, intersecting as you range over all $x\in X_{r_\phi}(\mathbb K)$, we see that $H':=\{h\in G: h(X_{r_\phi})= X_{r_\phi}\}$ is a Zariski closed subgroup. Now, $\phi$ acts on $X_{r_\phi}$ as multiplication by $g$, so that $g\in H'$. It follows that $H\subseteq H'$, and $X_{r_\phi}$ is $H$-invariant.

Let $\mathbb K$ be an uncountable algebraically closed field extending $F$.
As in the proof of Proposition~\ref{families}, it is still the case that $X_{r_\phi}$ has no nonconstant invariant rational functions over $\mathbb K$.
So, by~\cite{amerik-campana, BGR} as before, we have that~$\phi$, and hence the action of $H$ on $X$, admits a Zariski dense orbit over $\mathbb K$.
That is, $H$ acts generically transitively on $X_{r_\phi}$ over $\mathbb K$.

Now, consider $Z\subseteq H\times X_{r_\phi}^2$ given by
$$Z=\{(h,x,hx):h\in H, x\in X_{r_\phi}\}.$$
This is a family of invariant subvarieties of $(X_{r_\phi}\times X_{r_\phi},\phi\times\phi)$ because $H$ is commutative and $\phi$ is multiplication by $g\in H$.
It is a covering family because $H$ acts generically transitively on~$X_{r_\phi}$.
The members of this family are the graphs of the action of elements $h\in H$ on $X_{r_\phi}$, and so are of dimension $\dim X_{r_\phi}$.
Hence Proposition~\ref{families} gives us that the algebraic dimension of $(X_{r_\phi}\times X_{r_\phi},\phi_f\times\phi_f)$ is at least $\dim X_{r_\phi}$, as desired.
\end{proof}

\medskip
\subsection{Families of rational transformations} 
For this section we fix an algebraically closed field $k$.

Our main tool for recognising translations will be Proposition~\ref{findtv}, below, which asserts that it suffices to verify that the iterates of $\phi$ live in an algebraic family of rational transformations.
We leave our notion of families of rational transformations implicit in the statements of the lemmas and proposition of this section, but it is a natural notion modeled on the families of birational maps introduced by Demazure~\cite{demazure} (see also~\cite{blanc-furter}).

We will use the following semicontinuity lemma about algebraic families of rational transformations.
This is well known, though maybe somewhat delicate in positive characteristic.
We include proofs as we found no appropriate reference.

\begin{lemma}
\label{boundeg}
Suppose $X,Y,Z$ are irreducible varieties over an algebraically closed field~$k$, with $X\subseteq Z$ a subvariety.
Suppose $\psi:Y\times X\dto Z$ is a rational map, and let~$D$ be the set of points $a\in Y(k)$ such that $\psi_a:=\psi(a,-)$ is a rational transformation of~$X$.
Then there exists a positive integer $N$ such that the degree of the algebraic field extension $k(X)$ over $\psi_a^*k(X)$ is at most $N$, for all $a\in D$.
\end{lemma}

\begin{proof}
We will give an elementary argument that uses, mildly, some basic model theory of algebraically closed fields.

First of all, observe that $D$ is Zariski constructible.
Indeed, one obtains $D$ by first projecting onto $Y$ the domain of $\psi$ and then considering those $a$ such that $\psi_a:X\dto Z$ has its image in $X$, and finally restricting further to those $a$ such that $\psi_a:X\dto X$ is dominant.
All of these conditions are definable in the structure $(k,0,1,+,-,\times)$ and hence by quantifier elimination for the theory of algebraically closed fields (namely, Chevalley's theorem that the image of a Zariski constructible set is Zariski constructible) gives that $D$ is Zariski constructible.

We can therefore write
$$D=A_1\setminus B_1\cup\cdots\cup A_m\setminus B_m$$
where each $A_i\subseteq Y$ is an irreducible subvariety and $B_i\subsetneq A_i$ is a proper subvariety.
Let $U_i:=A_i\setminus B_i$, and consider the rational map $\psi^i:U_i\times X\dto Z$ obtained by restricting $\psi$ to $U_i$.
Note that, in fact, we have $\psi^i:U_i\times X\dto X$.
It suffices, of course, to prove the result for each $\psi^i$.
Hence, we may as well assume that $m=1$ and $D=Y(k)$ and $Z=X$.
That is, we have reduced to the case when $\psi:Y\times X\dto X$ and $\psi_a$ is a rational transformation of $X$ for all $a\in Y(k)$.

Note that as $k$ is algebraically closed, this remains true over any extension field $K\supseteq k$.
That is, for all $a\in Y(K)$, $\psi_a$ is a rational transformation of $X$ over~$K$.

Fix an algebraically closed field extension, $k\subseteq \KK$, that is large in the sense of having transcendence degree over $k$ greater than the cardinality of $k$.

Fix generators $f_1,\dots,f_d\in k(X)$.
Fix also one of the generators  $f:=f_i$ for some $i=1,\dots,d$.
Now fix $a\in Y(\KK)$.
Note that each $f_i\psi_a$ is a rational function on~$X$ over $k(a)$, and that $f_1\psi_a,\dots,f_d\psi_a$ generate $\psi_a^*\big(k(a)(X)\big)$.
As $\psi_a$ is dominant, this is a subfield of $k(a)(X)$ of the same transcendence degree over $k(a)$.
It follows that $f$ is algebraic over $\psi_a^*\big(k(a)(X)\big)$.
We obtain a polynomial $P^a(t,y_0,\dots y_d)\in k(t)[y_0,\dots,y_d]$ such that
\begin{itemize}
\item
$P^a(a,f,f_1\psi_a,\dots,f_d\psi_a)=0$ in the rational function field of $X$ over $k(a)$,
\item
$P^a(a,y_0,f_1\psi_a,\dots,f_d\psi_a)\neq 0$ as a polynomial in $y_0$ over the field
$$k(a,f_1\psi_a,\dots,f_d\psi_s)=\psi_a^*\big(k(a)(X)\big).$$
\end{itemize}
Consider the set $U^a\subseteq Y(\KK)$ of all $a'$ such that
\begin{itemize}
\item
$P^a(a',y_0,\dots,y_d)\in k(a')[y_0,\dots, y_d]$,
\item
$P^a(a',f,f_1\psi_{a'},\dots,f_d\psi_{a'})=0$ in $k(a')(X)$,
\item
$P^a(a',y_0,f_1\psi_{a'},\dots,f_d\psi_{a'})\neq 0$ in
$\psi_{a'}^*\big(k(a')(X)\big)[y_0]$.
\end{itemize}
Note that $U^a$ is a Zariski constructible set defined over $k$ that contains $a$.
Hence $\{U^a:a\in Y(\KK)\}$ is a Zariski constructible cover of $Y(\KK)$ defined over $k$.
By the large transcendence of $\KK$ over $k$ (namely {\em saturation} in model-theoretic terms) it follows that this cover has a finite subcover.
That is, $Y(\KK)=\bigcup_{i=1}^nU^{a_i}$ for some $a_1,\dots,a_n\in Y(\KK)$.
Set $N_i$ to be the maximum of the degrees in the variable $y_0$ of the polynomials $P^{a_1},\dots,P^{a_n}$.
(Remember that we had fixed $i$ and $f:=f_i$.)
Letting $N$ be the maximum of $N_1,\dots, N_d$, we get that  $[k(a)(X):\psi_a^*k(a)(X)]\leq N$, for all $a\in Y(\KK)$.
Applying this to all $a\in Y(k)$ proves the Lemma.
\end{proof}

Given a rational transformation $\phi:X\dto X$, by the {\em graph} of $\phi$ we will mean the (reduced) subvariety of $X^2$ obtained by taking the Zariski closure in $X^2$ of the set-theoretic graph of $\phi$ restricted to some nonempty Zariski open subset of $X$ on which it is defined.
This does  not depend on the Zariski open set chosen, and a rational map both uniquely determines and is determined by its graph.
The next lemma says that  the graphs of an algebraic family of transformations fit into finitely many families of subvarieties of $X^2$.
Let us first make precise our notion of families of subvarieties and fibres:

\begin{definition}
Suppose $X$ is an algebraic variety over an algebraically closed field~$k$.
By a {\em family of subvarieties of $X$} we mean an irreducible variety $V$ over~$k$ and an irreducible subvariety $\Gamma\subseteq V\times X$.
Given $a\in V(k)$, by the {\em fibre} $\Gamma_a$ we mean the set-theoretic (rather than scheme-theoretic) fibre; namely the reduced subvariety of $X$ such that $\Gamma_a(k)=\{b\in X(k):(a,b)\in \Gamma(k)\}$.
\end{definition}

\begin{lemma}
\label{boundgraph}
Suppose $X,Y,Z$ are irreducible varieties over an algebraically closed field~$k$, with $X\subseteq Z$ a subvariety.
Suppose $\psi:Y\times X\dto Z$ is a rational map, and let~$D$ be the set of tuples $a\in Y(k)$ such that $\psi_a:=\psi(a,-)$ is a rational transformation of~$X$.
Then there exist finitely many families of subvarieties of $X^2$ such that for all $a\in D$ the graph of $\psi_a$ appears as a fibre of one of these families.
\end{lemma}

\begin{proof}
Exactly as in the proof of Lemma~\ref{boundeg} we use the Zariski constructibility of~$D$ to reduce to the case that $D=Y(k)$ and $Z=X$.
That is, $\psi:Y\times X\dto X$ and $\psi_a$ is a rational transformation of $X$  for all $a\in Y(k)$.

Let $\Gamma(\psi)\subseteq Y\times X^2$ be the graph of $\psi$.
Note that $\Gamma(\psi)_a=\Gamma(\psi_a)$ for all points $a$ living in  some nonempty Zariski open subset $U\subseteq Y$.
This is because, for Zariski constructible sets, taking Zariski closures commutes with taking fibres, generically.
Hence, the graphs of $\psi_a$ for all $a\in U(k)$ live in a single family of subvarieties of $X^2$, namely the restriction of $\Gamma(\psi)$ to $U$.
We can now consider the restriction of $\psi$ to the irreducible components of $Y\setminus U$, and conclude by noetherian induction.
\end{proof}

Finally, let us record the fact that a family of subvarieties of a projective variety  involves only finitely many components of the Hilbert scheme.
Recall, from~\cite[page 265]{groth}, that for $X$ a projective variety over an algebraically closed field~$k$, $\hilb(X)$ is a scheme that parameterises the closed subschemes of $X$, and is made up of countably many irreducible components that are themselves projective varieties over~$k$.

\begin{lemma}
\label{boundhilb}
Suppose $\Gamma\subseteq V\times X$ is a family of subvarieties of a projective variety~$X$, over an algebraically closed field~$k$.
Then there exist finitely many components of $\hilb(X)$ such that, for all $a\in V(k)$, the fibre $\Gamma_a\subseteq X$ lives in one of these components.
\end{lemma}

\begin{proof}
By noetherian induction, it suffices to show that for some nonempty Zariski open subset $U\subseteq V$, the fibres $\Gamma_a$, for all $a\in U(k)$, land in the same component of $\hilb(X)$.
To do this, we use the fact that, after restricting to a Zariski open subset and taking a base extension, the co-ordinate projection $\Gamma\to V$ becomes flat with reduced scheme-theoretic fibres (so that there is no difference between the set-theoretic and scheme-theoretic fibre).
More precisely, there is a nonempty Zariski open subset $U\subseteq V$ and an irreducible affine variety~$\widetilde U$ with a surjective (purely inseparable) morphism $\pi:\widetilde U\to U$, such that if we look at the reduced fibre product,
$$\widetilde{\Gamma}:=\big(\Gamma\times_U\widetilde U\big)_{\operatorname{red}}\subseteq\widetilde U\times X$$
then the co-ordinate projection
$\widetilde{\Gamma}\to\widetilde U$ is flat with (geometrically) reduced scheme-theoretic fibres.
See~\cite[{Tag 0550}]{stacks-project}, for example.
The universal property of Hilbert schemes, together with irreducibility of $\widetilde U$, now gives rise to a morphism $\rho:\widetilde U\to W$, where $W$ is a component of $\hilb(X)$, such that $\rho(b)$ is $\widetilde{\Gamma}_b$ for each $b\in\widetilde U(k)$.
But $\widetilde{\Gamma}_b=\Gamma_{\pi(b)}$.
As~$\pi$ is surjective, we have that $\Gamma_a$ lives in $W$ for all $a\in U(k)$.
\end{proof}

With these preliminaries in place we can now present our main tool for finding translational dynamical systems.
The following Proposition is essentially an elaboration on and recasting of~\cite[Theorem~2.6]{cantat14}, and is based on Weil's regularisation (or ``group chunk") theorem from~\cite{weil}.

\begin{proposition}
\label{findtv}
Suppose $(X,\phi)$ is a rational dynamical system over an algebraically closed field~$k$ with the following property:
There exist $N\in\NN$, irreducible varieties~$Y$ and~$Z$ over~$k$, with $X\subseteq Z$, and a rational map $\psi:Y\times X\dto Z$ such that for each $r>N$ there is $a\in Y(k)$ with $\phi^r=\psi_a$.

Then $(X,\phi)$ is birationally equivalent to a translational dynamical system.
\end{proposition}

\begin{proof}
Passing to a projective closure, we may assume, without loss of generality, that $X$ is projective.

The first consequence of our assumption is that $\phi:X\dto X$ is in fact a birational transformation.
Indeed, by Lemma~\ref{boundeg}, the degree of the finite extension $k(X)/(\phi^r)^*k(X)$ remains bounded as $r>N$ varies.
But
$$[k(X):(\phi^r)^*k(X)]=r[k(X):\phi^*k(X)],$$
which forces $\phi^*k(X)=k(X)$.
That is, $\phi$ is birational.

Consider the group of birational transformations of $X$ generated by~$\phi$.
That is, let $\mathcal G:=\{\phi^r:r\in\mathbb Z\}$.
We view $\mathcal G$ as a set of $k$-points of $\hilb(X^2)$ by identifying each $\phi^r$ with its graph in $X^2$.
We claim that there are finitely many components of $\hilb(X^2)$ on which the elements of $\mathcal G$ live.
For $r$th iterates of $\phi$ where $r>N$ this follows directly from our assumption that they appear as fibres of $\psi$, by applying Lemmas~\ref{boundgraph} and~\ref{boundhilb}.
To deal with $\phi^{-r}$ for $r>N$, we consider the finitely many families of subvarieties of $X^2$ given by Lemma~\ref{boundgraph} applied to $\psi$, say
$$\Gamma_1\subseteq V_1\times X^2,\dots, \Gamma_n\subseteq V_n\times X^2.$$
For each $i=1,\dots, n$, let  $\Gamma_i^{\operatorname{opp}}\subseteq V_i\times X^2$ be the subvariety whose $k$ points are of the form $(v,x_1,x_2)$ where $(v,x_2,x_1)\in \Gamma_i$.
The graph of $\phi^{-r}$ for $r>N$ will appear as a fibre in one of these opposite families of subvarieties of $X^2$.
Applying Lemma~\ref{boundhilb} to these families, we therefore get that these $\phi^{-r}$ also involve only finitely many components of  $\hilb(X^2)$.
As this takes care of all but finitely many elements of~$\G$, we have now shown that $\G$ intersects only finitely many components of $\hilb(X^2)$ nontrivially.

At this point, at least when $k=\mathbb C$, that $(X,\phi)$ is translational follows formally from~\cite[Theorem~2.6]{cantat14}.
In fact the argument given there works over any algebraically closed field, and we repeat some of the details:
Let $\bir(X)$ be the group of birational transformations of $X$, viewed as a subset of $\hilb(X^2)$.
For each component $W$ of $\hilb(X^2)$, the intersection $\bir(X)\cap W$ is Zariski open in $W$.
Let $\overline{\mathcal G}$ be obtained by taking the Zariski closure of $\mathcal G$ in each of the finitely many $\bir(X)\cap W$ that arise.
The group law on $\bir(X)$ and its action as birational transformations on $X$ restrict to endow $\overline{\mathcal G}$ with the structure of a pre-group and $X$ with the structure of a pre-homogeneous space for~$\overline{\mathcal G}$.
Weil's regularisation theorem now produces a projective variety~$X'$, and a birational map $\psi:X'\dashrightarrow X$, over~$k$, such that $G:=\psi^{-1}\overline {\mathcal G}\psi$ is an algebraic subgroup of $\aut(X')$.
Letting $\phi':=\psi^{-1}\phi\psi\in G(k)$, therefore, we have that $(X',\phi')$ is translational and $\psi:(X',\phi')\dashrightarrow (X,\phi)$ is a birational equivalence of rational dynamical systems.
\end{proof}

The above method can also be used to show that translations are preserved under finite covers:

\begin{proposition}
\label{covers}
Suppose $(X,\phi)$ is a rational dynamical system over an algebraically closed field~$k$ that admits a generically finite-to-one equivariant dominant rational map $\pi:(X,\phi)\dashrightarrow(Y,\psi)$ where $(Y,\psi)$ is translational.
Then $(X,\phi)$ is birationally equivalent to a translational dynamical system.
\end{proposition}

\begin{proof}
Note, first of all, that $\phi$ must be a birational transformation of $X$.
Indeed, $k(X)$ is a finite algebraic extension of $k(Y)$, and $\phi^*$ restricts to $\psi^*$ on $k(Y)$.
As $(Y,\psi)$ is translational, $\psi^*$ is an automorphism.
If $\phi$ were not birational then the iterated images $(\phi^*)^nk(X)$ would form an infinite properly decreasing chain of subfields of $k(X)$ that contain $k(Y)$, which is impossible.

We may assume that $X$ is projective.

Consider the subgroup $\mathcal G=\{\phi^r:r\in\mathbb Z\}\leq\bir(X)$ generated by $\phi$.
As we saw in the proof of Proposition~\ref{findtv}, it suffices to show that $\mathcal G$ intersects only finitely many components of $\hilb(X^2)$ nontrivially.
As before, we are identifying a birational transformation with its graph in order to view $\bir(X)$ as living in $\hilb(X^2)$.
Note that, as $\pi$ is generically finite-to-one, dimension considerations imply that, for any $r\in\mathbb Z$, the graph of $\phi^r$, $\Gamma(\phi^r)$, is an irreducible component of $\pi^{-1}\big(\Gamma(\psi^r)\big)$.
But, as $\psi$ comes from an algebraic group action on $Y$, we know that
$$\mathcal H=\{\pi^{-1}\big(\Gamma(\psi^r)):r\in\mathbb Z\}$$
is a finite union of families of subvarieties of $X^2$.
It follows that the set of irreducible components of members of~$\mathcal H$ also fit into finitely many families of subvarieties of~$X^2$, and hence, by 
Lemma~\ref{boundhilb}, they involve only finitely many components of $\hilb(X^2)$.
\end{proof}

\bigskip
\section{The proof}
\label{sect:ftv}

\noindent
Our goal is to show that stabilised algebraic dimension is bounded by translation dimension.
We begin by showing that it is bounded by dimension; in fact we give a slightly more precise formulation of this for future use:

\begin{lemma}
\label{lem:specialise}
Suppose $(X,\phi)$ and $(Y,\psi)$ are rational dynamical systems over an uncountable algebraically closed field $k$, and 
$\theta_1,\dots,\theta_n\in\fix\big(k(X\times Y),\phi^*\otimes\psi^*\big)$ are algebraically independent over $\fix(k(Y),\psi^*)$.
Suppose $U\subseteq Y$ is a countable intersection of nonempty Zariski open subsets of $Y$ over~$k$.
Then there is $a\in U(k)$ such that $\theta_1(x,a),\dots,\theta_n(x,a)\in k(X)$ are algebraically independent over $k$.
\end{lemma}

\begin{proof}
It is a well known fact that if $(K,\sigma)\supseteq (L,\sigma)$ is a difference field extension, then $\fix(K,\sigma)$ is linearly disjoint from $L$ over $\fix(L,\sigma)$.
Indeed, toward a contradiction, suppose $c_1,\dots,c_n\in\fix(K,\sigma)$ are linearly independent over $\fix(L,\sigma)$ but $\sum_{i=1}^n a_ic_i=0$ is a nontrivial $L$-linear dependence.
We may assume that $n$ is minimal.
Dividing by $a_1$, which by minimality is not~$0$, we may assume that $a_1=1$.
Now, applying $\sigma$ yields $\sum_{i=1}^n \sigma(a_i)c_i=0$, and taking the difference of these linear relations yields an $L$-linear dependence among $c_2,\dots,c_n$, which by minimality must be trivial.
That is, $\sigma(a_i)=a_i$ for all $i=1,\dots,n$, contradicting the assumption that $c_1,\dots,c_n$ are linearly independent over $\fix(L,\sigma)$.

Let $L:=k(Y)$.
Applying the above observation to the difference field extension $\big(k(X\times Y),\phi^*\otimes\psi^*\big)\supseteq (L,\psi^*)$, induced by the co-ordinate projection $\pi:X\times Y\to Y$, we have that  $\fix\big(k(X\times Y),\phi^*\otimes\psi^*\big)$ is linearly disjoint from $L$ over $\fix(L,\psi^*)$.
So, as $\theta_1,\dots,\theta_n\in \fix\big(k(X\times Y),\phi^*\otimes\psi^*\big)$ were chosen algebraically independent over $\fix(L,\psi^*)$, they must remain algebraically independent over $L$.
That is,
$$\theta:=(\theta_1,\dots,\theta_n, \pi):X\times Y\dashrightarrow \AA^n\times Y$$
is a dominant rational map.
Let $k_0\subseteq k$ be a countable subfield over which $X,Y, \theta_1,\dots,\theta_n$ are defined, and such that  $U$ is the countable intersection of Zariski open subsets of $Y$ defined over $k_0$.
Let $a\in Y(k)$ be generic over $k_0$.
This exists because $k$ is an uncountable algebraically closed field.
Then $a\in U(k)$, and if we specialise $\theta$ to $a$, above, we obtain the dominant rational map 
$\theta_a=(\theta_1(x,a),\dots,\theta_n(x,a)):X\dashrightarrow \AA^n$.
This says precisely that $\theta_1(x,a),\dots,\theta_n(x,a)$ are algebraically independent over $k$ as rational functions on $X$.
\end{proof}

\begin{corollary}
\label{dimxbound}
Suppose $(X,\phi)$ is a rational dynamical system over an algebraically closed field $k$.
Then $\eadim(X,\phi)\leq\dim X$.
\end{corollary}

\begin{proof}
Suppose $(Y,\psi)$ is a rational dynamical system over~$k$ with
$$n:=\adim(X\times Y,\phi\times\psi)-\adim(Y,\psi).$$
We need to show that $n\leq\dim X$.
We have that $\fix\big(k(X\times Y),\phi^*\otimes\psi^*\big)$ is of transcendence degree $n$ over $\fix(k(Y),\psi^*)$.
So, let
$$\theta_1,\dots,\theta_n\in\fix\big(k(X\times Y),\phi^*\otimes\psi^*\big)$$
be algebraically independent over $\fix(k(Y),\psi^*)$.
Taking base extensions, we may assume that $k$ is uncountable.
By Lemma~\ref{lem:specialise}, there is $a\in Y(k)$ such that $\theta_1(x,a),\dots,\theta_n(x,a)\in k(X)$ are algebraically independent over $k$.
\end{proof}

Next we establish some elementary field theory preliminaries.

\begin{lemma} \label{lem:*}
Let $k\subseteq F$ be fields and let $q_1,\ldots ,q_s\in F[t]$ be a sequence of $k$-linearly independent polynomials.
Then there exists $p_1,\ldots ,p_s\in F[t]$ with the same $k$-span as $q_1,\ldots ,q_s$ but now satisfying the following property:
\begin{itemize}
\item[($\star$)]
for each $n$, if $\Lambda_n:=(p_i:\deg p_i=n)$ denotes the subsequence of $(p_1,\dots,p_s)$ of polynomials of degree exactly~$n$, and if $\Lambda_n\neq\emptyset$, then the sequence  of leading coefficients of the polynomials in~$\Lambda_n$ are $k$-linearly independent.
\end{itemize}
\end{lemma}

\begin{proof}
For any $k$-linearly independent sequence of nonzero $h_1,\ldots,h_s\in F[t]$, we let
$\Lambda_n(h_1,\ldots ,h_s):=(h_i:\deg h_i=n)$ and we let $V_n(h_1,\ldots ,h_s)$ denote the $k$-span of the leading coefficients of $\Lambda_n(h_1,\ldots ,h_s)$.
We define
$$e_n(h_1,\ldots ,h_s) = |\Lambda_n(h_1,\ldots ,h_s)| - {\rm dim}_kV_n(h_1,\dots,h_s),$$ where we take this quantity to be zero when $\Lambda_n(h_1,\ldots ,h_s)$ is empty.
Note that $e_n(h_1,\ldots ,h_s)$ is always a nonnegative integer, and our goal is to find a sequence of polynomials $p_1(t),\ldots ,p_s(t)\in F[t]$ with the same $k$-span as $q_1(t),\ldots ,q_s(t)$ such that $e_n(p_1,\ldots  ,p_s)=0$ for all $n$.

Let $\mathcal{S}\subseteq\omega^\omega$ denote the set of sequences of nonnegative integers that are eventually zero, equipped with the partial order $(a_i)\prec (b_i)$ if there is some $j$ such that $a_j<b_j$ and $a_i=b_i$ for $i>j$.
Since sequences in $\mathcal{S}$ are eventually zero, each element of $\mathcal{S}$ has only finitely many other elements of $\mathcal{S}$ that are less than it.  

Now among all $k$-linearly independent $h_1,\ldots ,h_s\in F[t]$ whose $k$-span is the same as that of $q_1,\ldots ,q_s$, choose $(p_1,\ldots ,p_s)$ with 
$(e_n(p_1,\ldots ,p_s):n\geq 0)$ minimal with respect to $\prec$ in $\mathcal S$.
If $e_n(p_1,\ldots ,p_s)=0$ for all $n$ then we are done, so we may assume, toward a contradiction, that this is not the case, and we let $m\geq 0$ be greatest such that $e_m(p_1,\ldots ,p_s)>0$.
By definition this means that if $\Lambda_m(p_1,\ldots, p_s) = (p_{m_1},..,p_{m_r})$ then there is some non-trivial $k$-linear combination of the form
$$c_1 p_{m_1}+\cdots + c_r p_{m_r} = h,$$
with $\deg h<m$.
Let us assume, for notational convenience, that $c_1$ is nonzero.
Now, for each $i=1,\dots,s$, set $\tilde{p}_i:=p_i$ if $i\neq m_1$ and $\tilde{p}_{m_1} := h$.
Then $\tilde{p}_1,\ldots ,\tilde{p}_s$ and $p_1,\ldots ,p_s$ have the same $k$-span.
Also, for $n>m$, we haven't changed anything; that is, $\Lambda_n(\tilde{p}_1,\ldots ,\tilde{p}_s)=\Lambda_m(p_1,\ldots ,p_s)$.
Hence, we still have $e_n(\tilde{p}_1,\ldots ,\tilde{p}_s) = 0$ for $n>m$.
On the other hand, by construction,
$$|\Lambda_m(\tilde{p}_1,\ldots ,\tilde{p}_s)|=r-1<|\Lambda_m(p_1,\ldots ,p_s)|$$
while
$$V_m(\tilde{p}_1,\ldots ,\tilde{p}_s)=V_m(p_1,\ldots ,p_s)$$
as the leading coefficient of the polynomial in $\Lambda_m(p_1,\ldots ,p_s)$ that we removed (namely $p_{m_1}$) was already a $k$-linear combination of the leadiong coefficients of the other elements of $\Lambda_m(p_1,\dots,p_s)$.
Hence $e_m(\tilde{p}_1,\ldots ,\tilde{p}_s) < e_m(p_1,\ldots ,p_s)$.
This contradicts the $\prec$-minimal choice of $(p_1,\ldots ,p_s)$.
\end{proof}

\begin{lemma}
\label{lem:val}
Suppose $k\subseteq F$ is a field extension of positive transcendence degree, and $p_1,\dots,p_s\in F(t)$ are rational functions that are $k$-linearly independent.
Then there are infinitely many $c\in F$ such that $p_1,\dots,p_s$ are defined at $c$ and $p_1(c),\dots,p_s(c)$ are $k$-linearly independent.
\end{lemma}

\begin{proof}
Clearing denominators, and avoiding the finitely many poles, it suffices to prove the lemma when $p_1,\dots,p_s\in F[t]$ are polynomials.
We may also assume that $F$ is a finitely generated extension of $k$.

Note that if we prove the result for some $q_1,\dots,q_s\in F[t]$ with the same $k$-span as $p_1,\dots,p_s$, then we will have proved it for $p_1,\dots,p_s$.
Indeed, for any $c\in F$, we have that $p_1(c),\dots,p_s(c)$ is obtained from $q_1(c),\dots,q_s(c)$ by applying an invertible $s\times s$ matrix over~$k$, so that the former is $k$-linearly independent if the latter is.
Hence, applying Lemma~\ref{lem:*}, and arguing by contradiction, if we suppose the result fails then there is a counterexample $p_1,\dots,p_s\in F[t]$ satisfying condition~($\star$) of Lemma~\ref{lem:*}.

Now, as $F$ is a nonalgebraic function field extension of $k$, it admits a nontrivial rank one discrete valuation $\nu$ with the property that $\nu(k^*)=\{0\}$.
(See, for example, \cite[Chapter~9]{AM}.)
Pick a finite-dimensional $k$-vector subspace $W$ of $F$ that contains all the coefficients of the~$p_i$.
Then there is some $d$ such that for all nonzero $w\in W$, $-d<\nu(w)< d$.

As $p_1(c),\dots, p_s(c)$ are linearly dependent over $k$ for all but finitely many $c\in F$, we can find such a $c\in F$ with $\nu(c)<-2d$.
Reordering if necessary, there is $1\leq \ell\leq s$, and $\lambda_1,\dots,\lambda_\ell\in k\setminus\{0\}$ such that
$$\sum_{i=1}^\ell \lambda_i p_i(c) = 0.$$
Now, let $N$ denote the maximum degree of $p_1,\ldots ,p_\ell$, and reorder so that $p_1,\dots,p_r$ have degree $N$, for some $1\leq r\leq \ell$, and the rest have degree strictly less than $N$.
Write
$p_i(t)=a_i t^N + q_i(t)$
with $q_i(t)$ of degree less than $N$, for each $i=1,\dots,\ell$.
So $a_i=0$ for $i=r+1,\dots,\ell$, and we get
$$\sum_{i=1}^r \lambda_i a_i c^N =  - \sum_{i=1}^\ell \lambda_i q_i(c).$$
By condition~($\star$), $a_1,\dots,a_r$ are $k$-linearly independent.
It follows that the left-hand side is $\theta c^N$ for some nonzero $\theta\in W$, and hence has value strictly less than $N\nu(c)+d$.
On the other hand, the right-hand side is a finite sum of terms of the form $wc^j$ with $j<N$ and $w\in W$, and so the valuation of each term on the right-side is greater than $(N-1)\nu(c)-d$.
It follows that $N\nu(c)+d>N\nu(c)-\nu(c)-d$, which gives $\nu(c)>-2d$, a contradiction.
\end{proof}

\begin{lemma}
\label{lem:vw}
Suppose $k\subseteq F_0\subseteq F$ are field extensions, and $V,W\subseteq F$ are finite dimensional $k$-vector subspaces.
Let
$$\frac{V}{W}:=\left\{f\in F : f=\frac{v}{w}\text{ for some }v\in V, 0\neq w\in W\right\}.$$
Then there are finite dimensional $k$-vector subspaces $V_0,W_0\subseteq F_0$ such that
$$\frac{V}{W}\cap F_0\subseteq\frac{V_0}{W_0}.$$
\end{lemma}

\begin{proof}
We may assume that $F$ is finitely generated over $F_0$ as we need only consider the subfield generated over $F$ by a $k$-basis for $V$ and $W$.
Descending one step at a time, it suffices to consider separately the cases when $F/F_0$ is purely inseparably algebraic, separably algebraic, and purely transcendental.

Suppose $F$ is a purely inseparable algebraic extension of $F_0$, and we are in characteristic $p>0$.
Let $(v_1,\dots,v_r)$ be a $k$-basis for $V$, and $(w_1,\dots,w_s)$ a $k$-basis for $W$.
Choose $n>0$ big enough so that each $w_i^{p^n}\in F_0$.
Let
$W_0:=\span_k\left\{w_1^{p^n},\dots,w_s^{p^n}\right\}$
and 
$$V_0:=F_0\cap\span_k\left\{v_iw_1^{a_1}\cdots w_s^{a_s}:i=1,\dots, r, a_1+\cdots+a_s=p^n-1\right\}.$$
These are finite-dimensional $k$-subspaces of $F_0$.
Given $z\in\frac{V}{W}\cap F_0$, write $z=\frac{v}{w}$ where $v\in V$ and $0\neq w\in W$.
Then
$$z=\frac{v}{w}=\frac{vw^{p^{n-1}}}{w^{p^n}}$$
Note that $w^{p^n}\in W_0$.
To see that $vw^{p^{n-1}}\in V_0$, note that it is in $F_0$, as both $z$ and $w^{p^n}$ are, and that if we write $v$ and $w$ out in terms of the basis then $vw^{p^{n-1}}$ has the correct form.

Next consider the case when $F$ is a separable algebraic extension of $F_0$.
Extending $F$ if necessary, we may assume that $F$ is a Galois extension of $F_0$.
Let $\id=\sigma_1,\sigma_2,\dots,\sigma_\ell$ be the distinct elements of $\gal(F/F_0)$.
Let $(v_1,\dots,v_r)$ be a $k$-basis for $V$, and $(w_1,\dots,w_s)$ a $k$-basis for $W$.
Let
$$W_0:=F_0\cap\span_k\left\{\prod_{j=1}^\ell\sigma_j(w_{i_j}): 1\leq i_1,\dots,i_\ell\leq s\right\}$$
and 
$$V_0:=F_0\cap\span_k\left\{v_i\prod_{j=2}^\ell\sigma_j(w_{i_j}): 1\leq i\leq r,\  1\leq i_2,\dots,i_\ell \leq s\right\}.$$
These are finite-dimensional $k$-subspaces of $F_0$.
Given $z\in\frac{V}{W}\cap F_0$, write $z=\frac{v}{w}$ where $v\in V$ and $0\neq w\in W$.
Then
$$z=\frac{v}{w}=\frac{v\prod_{j=2}^\ell\sigma_j(w)}{\prod_{j=1}^\ell\sigma(w)}.$$
Note that $\prod_{j=1}^\ell\sigma(w)\in F_0$, and so, as $z\in F_0$, we also have that $v\prod_{j=2}^\ell\sigma_j(w)\in F_0$.
Moreover, both have the correct forms to be in $W_0$ and $V_0$, respectively, so that $z\in\frac{V_0}{W_0}$, as desired.

Finally, suppose $F$ is a purely transcendental extension of $F_0$.
By induction, we may assume that $F$ is the field of rational function over $F_0$ in a single variable.
Again, let $(v_1,\dots,v_r)$ be a $k$-basis for $V$, and $(w_1,\dots,w_s)$ a $k$-basis for $W$.
Using Lemma~\ref{lem:val}, there exists $c\in F_0$ such that $v_1,\dots,v_r,w_1,\dots,w_s$ are all defined at $c$ and $w_1(c),\dots,w_s(c)$ remain $k$-linearly independent.
Let
$$V_0:=\span_k\{v_1(c),\dots,v_r(c)\}$$
and
$$W_0:=\span_k\{w_1(c),\dots,w_s(c)\}.$$
Given $d\in\frac{V}{W}\cap F_0$, write $d=\frac{v}{w}$ where $v\in V$ and $0\neq w\in W$.
Hence $dw=v$ as rational functions over $F_0$, and we can evaluate at $c$ to obtain $dw(c)=v(c)$.
Note that $v(c)\in V_0$, $w(c)\in W_0$, and $w(c)\neq 0$ as we chose $c$ such that $w_1(c),\dots,w_s(c)$ are $k$-linearly independent.
Hence $d=\frac{v(c)}{w(c)}\in \frac{V_0}{W_0}$.
\end{proof}

Now we are ready to prove the main theorem.

\begin{theorem}
\label{thm:tve}
Suppose $(X,\phi)$ is a rational dynamical system over an algebraically closed field $k$.
Then $\eadim(X,\phi)=\tdim(X,\phi)$.
\end{theorem}

\begin{proof}
We first show that $\eadim(X,\phi)\geq\tdim(X,\phi)$.
Let $\tdim(X,\phi)=n$ be witnessed by a translational dynamical system $(Y,\psi)$ over~$k$, of dimension~$n$, such that there is a dominant equivariant rational map $(X,\phi)\dashrightarrow (Y,\psi)$.
It follows that
\begin{eqnarray*}
\eadim(X,\phi)
&\geq&
\eadim(Y,\psi)\\
&=&\dim Y\ \ \text{ by Proposition~\ref{prop:sdim=dim} and Corollary~\ref{dimxbound}}\\
&=&
n
\end{eqnarray*}
as desired.

Now we work toward showing that $\eadim(X,\phi)\leq\tdim(X,\phi)$.
Our first step is to pass to an uncountable algebraically closed field extending~$k$.
In order to do so, let us first remark that, as $k$ is algebraically closed, both algebraic dimension and translation dimension are preserved under extension of the base field.
Indeed, both are given as the maximum dimension of a rational image of $X$ over~$k$ with certain properties -- 
algebraic dimension is the maximum dimension of an invariant fibration, and translation dimension is the maximum dimension of a translational dynamical system arising as a dominant equivariant image.
Hence, extending the base field certainly cannot drop these dimensions.
Morever, the witnesses to both algebraic and translation dimension are first-order definable properties of the parameters (in the language of rings), and hence a witness can be found over the algebraically closed field~$k$.
For example, if $f:(X,\phi)\dashrightarrow (Y,\id)$ is an invariant fibration over an algebraically closed field extension $K\supseteq k$, with $\dim Y\geq n$, then these facts, that $f$ is an invariant fibration and $\dim Y\geq n$, are expressed by a formula satisfied by the tuple of parameters used to define~$f$ and~$Y$, and so any realisation of that formula in~$k$ (and such exist as $k$ is algebraically closed) would witness that the algebraic dimension of $(X,\phi)$ is already at least~$n$ over~$k$.
Similarly, translation dimension over $K$ is already witnessed over $k$.
Now, stabilised algebraic dimension is the maximum of $\adim(X\times Y,\phi\times\psi)-\adim(Y,\psi)$ as $(Y,\psi)$ varies, and hence, as algebraic dimension is preserved, stabilised algebraic dimension does not drop if we extend the base field.
Hence, since translation dimension is also preserved, if we show that $\eadim(X,\phi)\leq\tdim(X,\phi)$ over $K\supseteq k$, then the inequality holds {\em a fortiori} over~$k$.
All this is to say that, for the purposes of proving this inequality, we can take a suitable base extension and assume that~$k$ is uncountable.

Now, suppose $(Y,\psi)$ is a rational dynamical system over~$k$ and
$$n:=\adim(X\times Y,\phi\times\psi)-\adim(Y,\psi).$$
We need to show that $(X,\phi)$ admits a dominant equivariant rational map to a translational dynamical system of dimension~$\geq n$.

Let $F:=k(X)$, $\sigma:=\phi^*$,  $L:=k(Y)$, and $\tau:=\psi^*$.
So we have $k$-difference fields $(F,\sigma)$ and $(L,\tau)$, and $k(X\times Y)=\Frac(F\otimes_kL)$ equipped with $\sigma\otimes\tau$.

As a matter of notation, given $k$-subspaces $V$ and $W$ of $F$ we let 
$$\frac{V}{W}:=\left\{f\in F : f=\frac{v}{w}\text{ for some }v\in V, 0\neq w\in W\right\}.$$

Fix $\theta_1,\dots,\theta_n\in\fix\big(k(X\times Y),\sigma\otimes\tau\big)$ algebraically independent over $\fix(k(Y),\tau)$.
Fix finite dimensional $k$-subspaces $V$ and $W$ of $F$ such that each $\theta_1,\dots,\theta_n$ can be written in the form
$\displaystyle\frac{\sum_{i=1}^rv_i\otimes g_i}{\sum_{j=1}^sw_j\otimes h_j}$ for some $v_1,\dots,v_r\in V$, $w_1,\dots,w_s\in W$, and $g_1,\dots,g_r,h_1,\dots,h_s\in L$.
Let
$$\Theta:=\{f\in F\colon f\in\frac{\sigma^m(V)}{\sigma^m(W)}\text{ for all sufficiently large }m\in\NN\}.$$

\begin{claim}
\label{claimz}
There exists a countable intersection of nonempty Zariski open subsets, $U\subseteq Y$, such that $\theta_\ell(x,a)\in \Theta$ for all $a\in U(k)$ and all $\ell=1,\dots,n$.
\end{claim}

\begin{proof}[Proof of Claim~\ref{claimz}]
It suffices to prove this for each fixed $\ell=1,\dots,n$.
Write
$$\theta_\ell = \frac{\sum_{i=1}^rv_i\otimes g_i}{\sum_{j=1}^sw_j\otimes h_j}$$
for some $v_1,\dots,v_r\in V$, $w_1,\dots,w_s\in W$, and $g_1,\dots,g_r,h_1,\dots,h_s\in L$.
Taking a minimal such presentation, so with $r$ and $s$ minimal, we may assume that the $w_1,\dots, w_s$ are linearly independent over $k$.
Fix $m\in\NN$.
Since $(\sigma\otimes\tau)^m(\theta_\ell)=\theta_\ell$, we have that
$$ \frac{\sum_{i=1}^rv_i\otimes g_i}{\sum_{j=1}^sw_j\otimes h_j}
=
\frac{\sum_{i=1}^r\sigma^m(v_i)\otimes\tau^m(g_i)}{\sum_{j=1}^s\sigma^m(w_j)\otimes\tau^m(h_j)}.$$
Let $a\in Y(k)$ be such that for all $m\in\NN$, $i=1,\dots,r$, and $j=1,\dots,s$, the rational functions $\tau^m(g_i)$ and $\tau^m(h_j)$ are defined at $a$, and it is not the case that $\tau^m(h_1)(a)=\cdots=\tau^m(h_s)(a)=0$.
This occurs off a countable union of proper Zariski closed sets.
Let $\lambda_{i,m}:=\tau^m(g_i)(a)$ and $\gamma_{j,m}:=\tau^m(h_j)(a)$.
We obtain:
$$ \theta_\ell(x,a)=\frac{\sum_{i=1}^r\lambda_{i,0}v_i}{\sum_{j=1}^s\gamma_{j,0}w_j}
=
\frac{\sum_{i=1}^r\lambda_{i,m}\sigma^m(v_i)}{\sum_{j=1}^s\gamma_{j,m}\sigma^m(w_j)}.$$
Here we are using linear independence of $w_1,\dots,w_s$, and the choice of $a$, to ensure that the above expressions are well-defined.

Noting that the $\lambda$'s and $\gamma$'s are in $k$ on which $\sigma$ acts trivially, we see that
$$ \theta_\ell(x,a)\in \frac{\sigma^m(V)}{\sigma^m(W)}.$$
This is for every $m\in\NN$.
Hence $\theta_\ell(x,a)\in \Theta$, as desired.
\end{proof}

Let $F_0:=k(\Theta)\subseteq F$.
By definition of $\Theta$, $F_0\subseteq F':= \bigcap_{m\ge 1} \sigma^m(F)$, and $\sigma$ restricts to a $k$-automorphism of $F'$.
We claim that $\sigma$ also restricts to an automorphism of $F_0$.
To see this, first observe that $\sigma(\Theta)\subseteq\Theta$ by definition.
On the other hand, if $g\in \Theta$ then $g\in F'$, and so we can consider $g':=\sigma^{-1}(g)\in F'$.
And one sees from the definition that $g'\in \Theta$.
So $\sigma$ permutes $\Theta$ and hence restricts to an automorphism of $F_0$.
%Additionally, by definition, for $g\in \Theta$ we have $\sigma^{-m}(g)\in (V/W)\cap F_0$ for all sufficiently large $m$.  

Combining Claim~\ref{claimz} with Lemma~\ref{lem:specialise}, we have that $F_0$ is a function field of transcendence degree at least $n$ over $k$.
Fix generators $f_1,\dots,f_d\in\Theta$ for $F_0$, and let $X_0\subseteq\mathbb A^d$ be an irreducible affine variety over $k$ such that $k[X_0]=k[f_1,\dots,f_d]$. 
We have a birational transformation $\phi_0:X_0\dashrightarrow X_0$ over~$k$ given by $\sigma|_{F_0}=\phi_0^*$, and the induced dominant rational map $X\dashrightarrow X_0$ is equivariant.
Moreover, $\dim X_0\geq n$.
It remains to show that $(X_0,\phi_0)$ is birationally equivalent to a translational dynamical system.
Actually, we will show the equivalent statement that $(X_0,\phi_0^{-1})$ is such.

We will use Proposition~\ref{findtv} to do this: we show that all sufficiently large iterates of $\phi_0^{-1}$ live in an algebraic family of rational transformations of~$X$.
Lemma~\ref{lem:vw} implies that 
$\frac{V}{W}\cap F_0\subseteq\frac{V_0}{W_0}$ for some finite dimensional $k$-subspaces $V_0, W_0$ of $F_0$.
Let $v_1,\dots,v_r$ be a $k$-linear basis for $V_0$ and $w_1,\dots,w_s$ a basis for $W_0$.
Consider the rational map
$\psi:\AA^{rd+sd}\times X_0\dto\AA^d$
given by 
$$\psi:=\left(\frac{\sum_{i=1}^ry_{i1}v_i}{\sum_{i=1}^sz_{i1}w_i},\dots,\frac{\sum_{i=1}^ry_{id}v_i}{\sum_{i=1}^sz_{id}w_i}\right).$$
Let $D$ be the set of all tuples $(\overline a,\overline b):=(a_{ij},b_{ij})\in k^{rd}\times k^{sd}$ such that $\psi_{(\overline a,\overline b)}$
is a rational transformation of $X_0$.
Note that we have $N\in\mathbb N$ such that for all $m>N$, 
$\phi_0^{-m}=(\sigma^{-m}(f_1),\dots,\sigma^{-m}(f_d))$
and each $\sigma^{-m}(f_i)\in \frac{V_0}{W_0}$.
That is, for each $m>N$, $\phi_0^{-m}=\psi_{(\overline a,\overline b)}$ for some $(\overline a,\overline b)\in D(k)$.
That is, $(X_0,\phi_0^{-1})$ satisfies the hypotheses of Proposition~\ref{findtv}, and hence is birationally equivalent to a translational dynamical system.
Hence so is $(X_0,\phi_0)$, as desired.
\end{proof}

\begin{corollary}
\label{maxsdim}
Suppose $(X,\phi)$ is a rational dynamical system over an algebraically closed field $k$.
Then $\eadim(X,\phi)=\dim X$ if and only if $(X,\phi)$ is birationally equivalent to a translational dynamical system.
\end{corollary}

\begin{proof}
The right-to-left implication is Proposition~\ref{prop:sdim=dim} together with Corollary~\ref{dimxbound}.
For the converse, suppose $\eadim(X,\phi)=\dim X$.
Then, by Theorem~\ref{thm:tve}, we have that $\tdim(X,\phi)=\dim X$, and hence $(X,\phi)$ is a generically finite cover of a translational dynamical system.
By Proposition~\ref{covers}, it follows that $(X,\phi)$ is itself birationally equivalent to a translational dynamical system.
\end{proof}

The following corollary was in fact our original motivation:

\begin{corollary}
\label{orthbound}
Suppose $(X,\phi)$ is a rational dynamical system over an algebraically closed field $k$.
Then $\eadim(X,\phi)>0$ if and only if $\adim(X^2,\phi\times\phi)>0$.

In particular, some cartesian power of $(X,\phi)$ has a nonconstant invariant rational function if and only if already the second power does.
\end{corollary}

\begin{proof}
Suppose $\adim(X^2,\phi\times\phi)>0$.
If $\adim(X,\phi)=0$ then, as $\eadim(X,\phi)\geq \adim(X^2,\phi\times\phi)-\adim(X,\phi)$, we get $\eadim(X,\phi)>0$.
On the other hand, if $\adim(X,\phi)>0$ then, as $\eadim(X,\phi)\geq\adim(X,\phi)$, we also get $\eadim(X,\phi)>0$.
This proves the right-to-left implication.

For the converse, suppose $\eadim(X,\phi)>0$.
Then, by Theorem~\ref{thm:tve}, we have a dominant equivariant rational map,  $$(X,\phi)\dashrightarrow (X_1,\phi_1),$$
with $(X_1,\phi_1)$ positive dimensional and translational.
If $\adim(X_1,\phi_1)>0$ then $\adim(X,\phi)>0$, and hence also $\adim(X^2,\phi\times\phi)>0$, as desired.
So we may assume $\adim(X_1,\phi_1)=0$.
Being translational, it follows from Proposition~\ref{translation}, that
$\adim(X^2_1,\phi_1\times\phi_1)\geq \dim X_1>0$.
But we have a dominant equivariant map $(X^2,\phi\times\phi)\dashrightarrow (X^2_1,\phi_1\times\phi_1)$, which will then witness that $\adim(X^2,\phi\times\phi)>0$.
\end{proof}

\vfill
\pagebreak
%\bibliographystyle{plain}
%\bibliography{dynamics}

\end{document}